\documentclass[12pt]{article}
\usepackage{amsfonts}
\usepackage{color}
\begin{document}

\newtheorem{theorem}{Theorem}[section]
\newtheorem{prop}[theorem]{Proposition}
\newtheorem{cor}[theorem]{Corollary}
\newtheorem{lem}[theorem]{Lemma}
\newtheorem{definition}[theorem]{Definition}
\newtheorem{ex}[theorem]{Example}
\newtheorem{no}[theorem]{Note}
\newtheorem{unnumber}{}
\renewcommand{\theunnumber}{\relax}
\newtheorem{prepf}[unnumber]{Proof}
\newenvironment{pf}{\prepf\rm}{\endprepf}
\newcommand{\qed}{\qquad$\star$}

\title{Counting false entries in truth tables of bracketed formulae connected by m-implication}
\author{Volkan Yildiz\\\\
\texttt{ah05146@qmul.ac.uk}\\
or\\
\texttt{vo1kan@hotmail.co.uk}
}
\date{}
\maketitle

\begin{abstract}
In this paper we count the number of rows $y_n$ with the value ``false'' in
the truth tables of all bracketed formulae with $n$ distinct variables connected by
the binary connective of ``modified-implication''. We find a recurrence and an asymptotic
formulae for $y_n$. We also determine the parity of $y_n$.
\end{abstract}

{\footnotesize

}

\section{Introduction}
In this paper we study enumerative and asymptotic questions on formulae
of propositional logic which are correctly bracketed chains of m-implications, 
where the letter `m' stands for `modified'. 

For brevity, we represent truth values of propositional variables 
and formulae by 1 for ``true'' and ``0'' for ``false''.

For background information on propositional logic the reader can  
refer to the following books,~\cite{DM}, and \cite{PJCS}, 
or to the introduction page of, \cite{P}.
In-fact this paper is an extension of  \cite{P}.
In \cite{P}, we have shown that the following results are true:

\begin{theorem}\label{t:f}
Let $f_n$ be the number of rows with the value ``false'' in the truth tables 
of all bracketed formulae with $n$ distinct propositions $p_1,\ldots,p_n$ 
connected by the binary connective of implication. Then 
\begin{equation}\label{e:f}
f_n =\sum_{i=1}^{n-1} (2^iC_i-f_i)f_{n-i}, \; \textit{ with } \; f_1=1
\end{equation}
and for large $n$, $\;f_n \sim \Bigg(\frac{3-\sqrt{3}}{6}\Bigg)\frac{2^{3n-2}}{\sqrt{\pi n^3}}$. 
Where $C_i$ is the $i$th Catalan number.
\end{theorem}
A number of new enumerative problems arise if we modify the 
binary connective of implication as in below cases.
\begin{itemize}
\item[Case(i)] Use $\rightharpoonup$ instead of $\to$, where  $\rightharpoonup$ defined as follows
\[
 \phi \rightharpoonup \psi \equiv  \phi \to \neg\psi 
\]
For any valuation $\nu$,
\[\nu(\phi\rightharpoonup\psi)=\cases{0 & if $\nu(\phi)=1$ and $\nu(\psi)=1$,\cr
1 & otherwise.\cr}\]
\item[Case(ii)] Use $\leftharpoonup$ instead of $\to$, where  $\leftharpoonup$ defined as follows
\[
 \phi \leftharpoonup \psi \equiv  \neg \phi \to\psi 
\]
For any valuation $\nu$,
\[\nu(\phi\leftharpoonup\psi)=\cases{0 & if $\nu(\phi)=0$ and $\nu(\psi)=0$,\cr
1 & otherwise.\cr}\]
\item[Case(iii)] Use $\rightleftharpoons$ instead of $\to$, where  $\rightleftharpoons$ defined as follows
\[
 \phi \rightleftharpoons \psi \equiv  \neg \phi \to \neg\psi 
\]
For any valuation $\nu$,
\[\nu(\phi\rightleftharpoons \psi)=\cases{0 & if $\nu(\phi)=0$ and $\nu(\psi)=1$,\cr
1 & otherwise.\cr}\]
\end{itemize}
Let  $s_n$, $h_n$ be the number of rows with the value ``false'' in the truth tables 
of all bracketed formulae with $n$ distinct propositions $p_1,\ldots,p_n$ 
connected by the binary connective of m-implication, in the case (iii) and (ii), respectively.

\subsection{Case(iii)}
A row with the value false comes from an expression $\psi\rightleftharpoons \chi$ where $\nu(\psi)=0$
and $\nu(\chi)=1$. If $\psi$ contains $i$ variables, then $\chi$ contains
$n-i$, and the number of choices is given by the summand:
\begin{equation}\label{e:s}
s_n = \sum_{i=1}^{n-1} s_i (2^{n-i}C_{n-i}-s_{n-i}), \textit{ where } s_0=0, s_1=1.
\end{equation}
The recurrence relation (\ref{e:s}) is equivalent to the recurrence relation (\ref{e:f}),
 so all the results we have in \cite{P}, and \cite{PP} hold for the case(iii) too. 

\subsection{Case(ii)}
 A row with the value false comes from an expression $\psi\leftharpoonup \chi$ where $\nu(\psi)=0$
and $\nu(\chi)=0$. If $\psi$ contains $i$ variables, then $\chi$ contains
$n-i$, and the number of choices is given by the summand:
\begin{equation}\label{e:h}
h_n=\sum_{i=1}^{n-1} h_ih_{n-i}, \textit{ where } h_0=0, h_1=1.
\end{equation}
The recurrence relation (\ref{e:h}) is very well known; it is the recurrence relation for Catalan numbers. 
\begin{cor}
Suppose  we have all possible well-formed formulae obtained from
$p_1 \leftharpoonup p_2 \leftharpoonup \ldots \leftharpoonup p_n$ by inserting brackets, 
where $p_1,\ldots,p_n$ are distinct propositions. 
Then each formula defines the same truth table.
\end{cor}

\begin{ex}
 Here are the truth tables, (merged into one), for the bracketed m-implications, in $n=3$ variables.
\[
\begin{array}{|l|l|l|c|c|}
\hline p_1 & p_2 & p_3 & p_1\leftharpoonup (p_2 \leftharpoonup p_3) & (p_1 \leftharpoonup p_2)\leftharpoonup p_3 \\
\hline 1 & 1 & 1 & 1 & 1 \\
\hline 1 & 1 & 0 & 1 & 1 \\
\hline 1 & 0 & 1 & 1 & 1 \\
\hline 1 & 0 & 0 & 1 & 1 \\
\hline 0 & 1 & 1 & 1 & 1 \\
\hline 0 & 1 & 0 & 1 & 1 \\
\hline 0 & 0 & 1 & 1 & 1 \\
\hline 0 & 0 & 0 & {\color{blue} 0} & {\color{blue} 0} \\
\hline
\end{array}
\]
\end{ex}

\subsection{Case(i)} We are interested in 
{\em bracketed m-implications, case(i)}, which are formulae obtained from 
$p_1\rightharpoonup p_2 \rightharpoonup \ldots \rightharpoonup p_n$ 
by inserting brackets so that the result is well-formed, 
where $p_1,\ldots ,p_n$ are distinct propositions.

\begin{prop}\label{p:1}
Let $y_n$ be the number of rows with the value ``false'' in the truth tables of all brackted 
m-implications, case(i), with $n$ distinct variables. Then 
\begin{equation}\label{e:y}
y_n=\sum_{i=1}^{n-1} \bigg((2^iC_i - y_i)(2^{n-i}C_{n-i}-y_{n-i})\bigg), \textit{ with } y_0=0, \; y_1=1.
\end{equation}
\end{prop}
\begin{pf}
A row with the value false comes from an expression $\phi\rightharpoonup \psi$, 
where $\nu(\phi)=1$ and $\nu(\psi)=1$. If $\phi$ contains $i$ variables, then 
$\psi$ contains $n-i$ variables, and the number of choices is given by the summand in 
the proposition.
\qed
\end{pf}

\begin{ex}\label{e:q}
\[
y_1= 1, y_2=(2^1C_1-y_1)(2^1C_1-y_1)=1
\]
and 
\[
y_3=(2^1C_1-y_1)(2^2C_2-y_2) + (2^2C_2-y_2)(2^1C_1-y_1) = 3+3= 6.
\]
\end{ex}

\begin{ex} \label{e:1}
Here are the truth tables, (merged into one), for the two bracketed m-implications, case(i),
in $n=3$ variables. Where the corresponding rows with the value false are in blue:
\[
\begin{array}{|l|l|l|c|c|}
\hline p_1 & p_2 & p_3 & p_1\rightharpoonup (p_2 \rightharpoonup p_3) & (p_1 \rightharpoonup p_2)\rightharpoonup p_3 \\
\hline 1 & 1 & 1 & 1 & 1 \\
\hline 1 & 1 & 0 & 1 & {\color{blue} 0} \\
\hline 1 & 0 & 1 & {\color{blue} 0} & {\color{blue} 0} \\
\hline 1 & 0 & 0 & 1 & {\color{blue} 0} \\
\hline 0 & 1 & 1 & {\color{blue} 0} & 1 \\
\hline 0 & 1 & 0 & 1 & 1 \\
\hline 0 & 0 & 1 & {\color{blue} 0} & 1 \\
\hline 0 & 0 & 0 & 1 & 1 \\
\hline
\end{array}
\]
which coincides with the result we had from Example~\ref{e:q}.
\label{e:t}
\end{ex}

Using Proposition~\ref{p:1}, it is straightforward to calculate the values of $y_n$ 
for small $n$. The first $22$ values are
\begin{eqnarray*}
\{y_n\}_{n\geq 1} &=& 1, 1, 6, 29, 162, 978, 6156, 40061, 267338, 1819238,\\ 
&& 12576692, 88079378, 623581332, 4455663876, 32090099352, \\
&& 232711721757, 1697799727066, 12452943237342, 91774314536100,\\
&&  679234371006982,5046438870909244, 37623611703611452, \ldots
\end{eqnarray*}

\section{Generating Function}
Recall from \cite{P}, that the number of bracketings of a product of $n$ terms
is the Catalan number with the generating function
\[C_n= \frac{1}{n}{2n-2\choose n-1}, \textit{ with } C_0=0, \;\; \sum_{n\ge1}C_nx^n = (1-\sqrt{1-4x})/2 \]
respectively (see also~\cite[page 61]{PJCC}).

Let $g_n$ be the total number of rows in all truth tables 
for bracketed m-implications, case(i), with $n$ distinct variables. 
It is clear that $g_n=2^nC_n$, with $g_0=0$.
Let $Y(x)$ and $G(x)$ be the generating functions for $y_n$, and $g_n$, respectively. 
That is, $Y(x)= \sum_{n\geq 1} y_nx^n$, and $G(x)=\sum_{n\geq 1} g_nx^n$ .\\\\
Since,
\[ y_n=\sum_{i=1}^{n-1} \bigg((2^iC_i - y_i)(2^{n-i}C_{n-i}-y_{n-i})\bigg), \;\;\textit{ where }\; y_0=0, \; y_1=1. \]
Then,{\small
\begin{eqnarray*}
\sum_{n\geq1}y_n x^n  &=&  x + \sum_{n\geq1}\sum_{i=1}^{n-1} 2^iC_i2^{n-i}C_{n-i} x^n -  \sum_{n\geq1}\sum_{i=1}^{n-1} 2^iC_i y_{n-i} x^n - \\ 
&& \sum_{n\geq1}\sum_{i=1}^{n-1}y_i 2^iC_{n-i}x^{n-i} +  \sum_{n\geq1}\sum_{i=1}^{n-1} y_{i}y_{n-i} x^n 
\end{eqnarray*}}
Now it is straightforward to get the following result:
\begin{equation}\label{eq:1}
Y(x) = x+(G(x)-Y(x))^2
\end{equation}
where $G(x)$ can be obtained from the generating function of $C_n$ by replacing
$x$ by $2x$: that is,
\begin{equation}\label{eq:2}
G(x) = (1- \sqrt{1-8x})/2.
\end{equation}
Substituting (\ref{eq:2}) into  (\ref{eq:1}) gives the following quadratic equation:
\begin{equation}\label{eq:3}
2Y(x)^2+2Y(x)(\sqrt{1-8x}-2)+(1-\sqrt{1-8x}-2x)=0
\end{equation}
Solving equation~(\ref{eq:3}) gives the following proposition:
\begin{prop}
The generating function for the sequence $\{y_n\}_{n\geq 1}$ is given by
\[
Y(x) = \frac{2-\sqrt{1-8x} - \sqrt{3-4x-2\sqrt{1-8x}}}{2}.
\]
\label{p:gf}
\end{prop}
(As with the Catalan numbers, the choice of sign in the square root is made to
ensure that $Y(0)=0$.)
With the help of Maple we can obtain the first $22$ terms of the above series, 
and hence give the first $22$ values of $y_n$; these agree with the values
found from the recurrence relation.

\section{Asymptotic Analysis}

In this section we want to get an asymptotic formula for the coefficients of
the generating function $Y(x)$ from Proposition~\ref{p:gf}. 
We use the following result~\cite[page 389]{FAC}:

\begin{prop}
Let $a_n$ be a sequence whose terms are positive for sufficiently large $n$.
Suppose that $A(x)=\sum_{n\geq 0} a_nx^n$ converges for some value of $x>0$.
Let $f(x)= (-\ln(1-x/r))^b(1-x/r)^c$, where $c$ is not a
positive integer, and we do not have $b=0$ and $c=0$. 
Suppose that $A(x)$ and $f(x)$ each have a singularity at $x=r$ 
and that $A(x)$ has no singularities in the interval $[-r,r)$. 
Suppose further that $\lim_{x\to r} \frac{A(x)}{f(x)}$ exists and has nonzero value $\gamma$. Then 
\[a_n \sim \cases{
\gamma{n-c-1\choose n} (\ln{n})^br^{-n}, & if $c\not=0$,\cr\cr
\frac{\gamma b(\ln{n})^{b-1}}{n}, & if $c=0$.\cr}\]
\label{p:b}
\end{prop}

\begin{no}\label{n:b}
We also have
\[
{n-c-1\choose n}\sim \frac{n^{-c-1}}{\Gamma(-c)},
\]
where the standard gamma-function
\[\Gamma(x)=\int_0^\infty t^{x-1}\mathrm{e}^{-t}\,\mathrm{d}t, \;\textit{ with } \;\Gamma(x+1)=x\Gamma(x), \; \Gamma(1/2)=\sqrt{\pi} .\]
It follows that $\Gamma(-1/2)=-\sqrt{\pi}/2$ .
\end{no}
Recall that $G(x)=(1-\sqrt{1-8x})/2$, therefore 
\[Y(x)= \frac{(1+2G(x))-\sqrt{(1+4G(x))-4x}}{2}.\]
As in \cite{P}, before studying $Y(x)$, we first study $G(x)$. This $G(x)$ could easily be
studied by using the explicit formula for its coefficients, which is
$2^n{2n-2\choose n-1}/n$. But our aim is to understand how to handle the
square root singularity. A square root singularity occurs while attempting
 to raise zero to a power which is not a positive integer. 
Clearly the square root, $\sqrt{1-8x}$, has a singularity at $1/8$.
Therefore by Proposition~\ref{p:b}, $r=1/8$. We have
$G(1/8)=1/2$, so we would not be able to divide $G(x)$ by a
suitable $f(x)$ as required in Proposition~\ref{p:b}. To create a function
which vanishes at $\frac{1}{8}$, we simply look at
$A(x)=G(x)-1/2$ instead. That is, let
\[f(x)=(1-x/r)^{1/2} =  (1-8x)^{1/2}.\]
Then
\[\gamma=\lim_{x\rightarrow1/8}\frac{A(x)}{\sqrt{1-8x}} = -\frac{1}{2}.\]

Now by using Proposition~\ref{p:b} and Note~\ref{n:b},
\[
g_n \sim -\frac{1}{2}{n-\frac{3}{2}\choose n} \bigg(\frac{1}{8}\bigg)^{-n}  \sim -\frac{1}{2}\,\frac{8^n n^{-3/2}}{\Gamma (-1/2)} = \frac{2^{3n-2}}{\sqrt{\pi n^3}}.
\]

We are now ready to tackle $Y(x)$, and state the main theorem of the paper.

\begin{theorem}
Let $y_n$ be number of rows with the value false in the truth tables of all the
bracketed m-implications, case(i), with $n$ distinct variables. Then
\[y_n \sim \left(\frac{10-2\sqrt{10}}{10}\right)\frac{2^{3n-2}}{\sqrt{\pi n^3}}.\]
\end{theorem}
\begin{pf}
Recall that
\[
Y(x) = \frac{2-\sqrt{1-8x} - \sqrt{3-4x-2\sqrt{1-8x}}}{2}.
\] 
We find that $r=\frac{1}{8}$, and $f(x)=\sqrt{1-8x}$. 
Since $Y(1/8)=(2\sqrt{2}-\sqrt{5})/2\sqrt{2} \not= 0$, we need a function 
which vanishes at $Y(1/8)$, thus we let $A(x)=Y(x)-Y(1/8)$.

\[\lim_{x\to1/8} \frac{A(x)}{f(x)} 
= \lim_{x\to 1/8}  \frac{-\sqrt{2}\sqrt{1-8x} - \sqrt{2}\sqrt{3-4x-2\sqrt{1-8x}}+\sqrt{5}}{2\sqrt{2}\sqrt{1-8x}} .\]
Let $v=\sqrt{1-8x}$. Then
\begin{eqnarray*}
\gamma &=& \lim_{v\to 0} \frac{-\sqrt{2}v - \sqrt{v^2-4v+5} + \sqrt{5}}{2\sqrt{2}v} = \lim_{v\to 0} \frac{-\sqrt{2} - \frac{1}{2}(2v-4)(v^2-4v+5)^{\frac{-1}{2}}}{2\sqrt{2}} \\
&=& \frac{-\sqrt{2}+\frac{2}{\sqrt{5}}}{2\sqrt{2}} =  -\frac{10-2\sqrt{10}}{20},
\end{eqnarray*}
where we have used l'H\^opital's Rule in the penultimate line.

Finally,
\[y_n \sim  -\frac{10-2\sqrt{10}}{20}{n-\frac{3}{2}\choose n}
\left(\frac{1}{8}\right)^{-n} \sim \left( \frac{10-2\sqrt{10}}{10}\right)
\frac{2^{3n-2}}{\sqrt{\pi n^3}},\]
and the proof is finished. \qed
\end{pf}

The importance of the constant $ \frac{10-2\sqrt{10}}{10} =0.367544468$ 
lies in the following fact:

\begin{cor}
Let $g_n$ be the total number of rows in all truth tables for bracketed
m-implications, case(i), with $n$ distinct variables, and $y_n$ the number of rows with the
value ``false''. Then $\lim_{n\to\infty}y_n/g_n= \frac{10-2\sqrt{10}}{10}$.
\end{cor}

The table below illustrates the convergence.

\[
\begin{array}{|c|c|c|c|}
\hline n & y_n & g_n & y_n/g_n \\
\hline 1 & 1 & 2 & 0.5 \\
\hline 2 & 1 & 4 & 0.25 \\
\hline 3 & 6 & 16 & 0.25 \\
\hline 4 & 29 & 80 & 0.3625 \\
\hline 5 & 162 & 448 & 0.36160714286\\
\hline 6 & 978 & 2688 & 0.36383928571\\
\hline 7 & 6156 & 16896 & 0.36434659091 \\
\hline 8 & 40061 & 109824 & 0.36477454837\\
\hline 9 & 267338 & 732160 & 0.36513603584 \\
\hline 10 & 1819238 & 4978688 & 0.36540510271 \\
\hline 100 & - &- & 0.36735248210 \\
\hline
\end{array}\]

\begin{cor}
Let \[ P(y_n)= \frac{y_n}{g_n} \; \textit{ and } \;  P(f_n)= \frac{f_n}{g_n}\]  
then we have the following inequality
\[
P(y_n) \geq P(f_n) .
\]
Where $f_n$ is defined in Theorem~\ref{t:f}.
\end{cor}

\begin{cor} \label{c:c}
Let $d_n$ be the number of rows with the value ``true'' in the truth tables of all 
bracketed formulae with $n$ distinct variables connected by the binary 
connective of m-implication, case(i). Then 
\[
d_n= g_n -y_n, \textit{ with } t_0 =0,
\]
and for large $n$, 
\[
d_n \sim \Bigg(\sqrt{\frac{2}{5}}\Bigg) \frac{2^{3n-2}}{\sqrt{\pi n^3}} .
\]
\end{cor}
Using this Corollary~\ref{c:c}, it is straightforward to calculate the values of $d_n$. 
The table below illustrates this up to $n=10$.
\[
\begin{array}{|l|c|c|c|c|c|c|c|c|c|c|c|}
\hline n&0 & 1& 2 & 3 & 4 & 5 & 6 & 7 & 8 & 9 & 10\\
\hline d_n & 0 & 1 & 3 & 10  & 51 & 286 & 1710 & 10740 & 69763 & 464822 & 3159450\\
\hline
\end{array}
\]

\section{Parity}
For brevity, we represent the set of even counting numbers by 
the capital letter $E$, the set of odd counting numbers by the capital letter $O$, 
and the set of natural numbers, $\{1,2,3,4,...\}$, by $\mathbb{N}$.
 
We begin by determining the parity of Catalan number $C_n$, which has the following recurrence relation
\begin{equation}\label{e:c}
C_n=\sum_{i=1}^{n-1}C_iC_{n-1}, \; \textit{ with } C_0=0, C_1=1.
\end{equation}

From the Segner's recurrence relation, $C_n$ can be expressed as a piecewise function, 
with respect to the parity of $n$, (see~\cite[page 329]{C}).

\[C_n=\cases{2(C_1C_{n-1}+C_2C_{n-2}+\ldots+C_{\frac{n-1}{2}}C_{\frac{n+1}{2}}) \;\;& if $n\in O$,\cr\cr
2(C_1C_{n-1}+C_2C_{n-2}+\ldots+C_{\frac{n-2}{2}}C_{\frac{n+2}{2}})+C_{\frac{n}{2}}^2 \;\; & if $n\in E$.\cr}\]

\begin{lem}[Parity of $C_n$]  ~\cite{PP} \label{l:1}
\[ C_n\in O \Longleftrightarrow n=2^i, \textit{ where } i\in\mathbb{N} . \]
\end{lem}
\begin{pf}
\[
\textit{ For } 
n\geq 2, \; C_n\in O \Longleftrightarrow C_{\frac{n}{2}}^2 \in O  \Longleftrightarrow C_{\frac{n}{2}}\in O \Longleftrightarrow n=2^i \;\;\forall i\in \mathbb{N}. 
\]
Note that $C_1=1\in O$. \qed 
\end{pf}

By using Proposition~\ref{p:1}, we get the following triangular table. 
Where the left hand side column represents the sum of the corresponding row.
\begin{center}
\begin{tabular}{rccccccccc}
$y_2$:&    &    &    &    &  1\\\noalign{\smallskip\smallskip}
$y_3$:&    &    &    &  3 &    & 3  \\\noalign{\smallskip\smallskip}
$y_4$:&    &    &  10 &    &  9 &    &  10 \\\noalign{\smallskip\smallskip}
$y_5$:&    &  51 &    &  30 &    &  30 &    &  51 \\\noalign{\smallskip\smallskip}
$y_6$:&  286 &    &  153 &    &  100 &    &  153 &    &  286\\\noalign{\smallskip\smallskip}
\end{tabular}
\end{center}
\begin{theorem}[Parity of $y_n$]\label{f:c} The sequence $\{y_n\}_{n\geq 1}$ preserves the parity of $C_n$.
\end{theorem}

\begin{pf} If an additive partition of $y_n$, (which is determined by the recurrence 
relation~(\ref{e:y})), is odd, then it comes as a pair; i.e.  
\[
(2^iC_i-f_i)(2^{n-i}C_{n-i}-y_{n-i}) \in O \Longleftrightarrow y_i, y_{n-i}.
\]
Hence, $\bigg((2^iC_i-y_i)(2^{n-i}C_{n-i}-y_{n-i})+(2^{n-i}C_{n-i}-y_{n-i})(2^iC_i-y_i)\bigg) \in E.$

Thus, $y_n$ can be expressed as a piecewise function depending on the parity of $n$:
\[y_n=\cases{2\sum_{i=1}^{\frac{n-1}{2}} ((2^iC_i-y_i)(2^{n-i}C_{n-i}-y_{n-i}) ) \;\;& if $n\in O$,\cr\cr
\bigg(2\sum_{i=1}^{\frac{n-2}{2}} ((2^iC_i-y_i)(2^{n-i}C_{n-i}-y_{n-i}) ) \bigg)+ (2^{\frac{n}{2}}C_{\frac{n}{2}}-y_{\frac{n}{2}})^2\;\; & if $n\in E$.\cr}\]
Finally,  
\[
y_n\in O \Longleftrightarrow (2^{\frac{n}{2}}C_{\frac{n}{2}}-y_{\frac{n}{2}})^2 \in O \Longleftrightarrow y_\frac{n}{2} \in O  \Longleftrightarrow  n=2^i, \;\; \forall i\in\mathbb{N}.
\]
Note that $y_1=1\in O$. \qed
\end{pf}

\begin{prop}[Parity of $d_n$]
The sequence $\{d_n\}_{n\geq 1}$ preserves the parity of $C_n$.
\end{prop}
\begin{pf}Since 
\[
d_n=g_n-y_n=2^n C_n-y_n, \textit{ with } n\geq 1 
\]
The sequence $\{g_n\}_{n\geq 1}$ is always even, and the sequence 
$\{y_n\}_{n\geq 1}$ preserves the  parity of $C_n$ by Theorem~\ref{f:c}. 
Therefore the sequence $\{d_n\}_{n\geq 1}$ preserves the parity of $C_n$.
\qed
\end{pf}

\section{A fruitful tree}
We begin by recalling following definitions:
\begin{definition}\cite{PP}, 
The nth {\bf Catalan tree}, $A_n$, is a combinatorical object, 
characterized by one root, $(n-1)$  main-branches, and $C_n$  sub-branches.  
Where each main-branch gives rise to a number of sub-branches, and 
the number of these sub-branches is determined by the additive 
partition of the corresponding Catalan number,
 as determined by the recurrence relation (\ref{e:c}).
\end{definition}

\begin{definition}\cite{PP}, 
The Catalan tree $A_n$ is {\bf fruitful} iff each sub-branch of $A_n$ has fruits. 
We denote this new tree by $A_n(\mu_i)$, 
where $\{\mu_i\}_{i\geq 1}$ is the corresponding fruit sequence.
\end{definition}

\begin{ex}

Let $\{y_n\}_{n\geq 1}$ be the corresponding fruit sequence for the Catalan tree $A_n$. 
Then $A_n(y_n)$ has the following symbolic representation, 
{\small
\begin{eqnarray*}
&\big((2^1C_1-f_1)(2^{n-1}C_{n-1}-y_{n-1}), \ldots, (2^{n-1}C_{n-1}-y_{n-1})(2^1C_1-f_1)\big) \\
& (C_1C_{n-1},C_2C_{n-2},\ldots,C_{n-2}C_2,C_{n-1}C_1) \\
&  (1,1,\ldots,1,1) \\
& (1).
\end{eqnarray*}
}
\end{ex}

\begin{ex}
Let $\{d_n\}_{n\geq 1}$ be the corresponding fruit sequence for the Catalan tree $A_n$. 
Then $A_n(d_n)$  has the following symbolic representation, 
{\small
\begin{eqnarray*}
&\big((2^n-(2^1C_1-f_1)(2^{n-1}C_{n-1}-y_{n-1})), \ldots, (2^n- (2^{n-1}C_{n-1}-y_{n-1})(2^1C_1-f_1) )\big) \\
& (C_1C_{n-1},C_2C_{n-2},\ldots,C_{n-2}C_2,C_{n-1}C_1) \\
&  (1,1,\ldots,1,1) \\
& (1).
\end{eqnarray*}
}
\end{ex}

\begin{prop}\label{p:s}  For $n>1$, let $a_n(y_n)$ and $a_n(d_n)$ be the total number of components of 
the fruitful trees $A_n(y_n)$ and $A_n(d_n)$ respectively. Then  
\[
a_n(y_n) = y_n + C_n + n, \; \textit{ and } \; a_n(d_n) = d_n + C_n + n .
\]
\end{prop}

Using Proposition~\ref{p:s}, it is straightforward to calculate the values of $a_n(y_n)$ , and $a_n(d_n)$. 
The table below illustrates this up to $n=10$.

\[
\begin{array}{|l|c|c|c|c|c|c|c|c|c|c|c|}
\hline n & 0 & 1& 2 & 3 & 4 & 5 & 6 & 7 & 8 & 9 & 10\\
\hline a_n(y_n)  & 0 & 2 & 4 & 11 & 38 & 181 & 1026 & 6295 & 40498 & 268777 & 1824110 \\
\hline a_n(d_n)  & 0 & 2 & 6 & 15 & 60 & 305 & 1758 & 10879 & 70200 & 466261 & 3164322 \\
\hline
\end{array}
\]

\begin{cor}For $n>1$,  $a_n(y_n)$, and $a_n(d_n)$  are odd iff $n\in O$. \end{cor}
\begin{pf}
Since,
\[
a_n=(C_n+n) \in O \Longleftrightarrow n\in O \textit{ or } n=2^i, \textit{ and } y_n,d_n \in O \Longleftrightarrow n= 2^i  .
\]
Therefore, $a_n(y_n), a_n(d_n)\in O \Longleftrightarrow n\in O$ .
\end{pf}

\begin{table}[b]
\begin{verse}
Onlar ki kurtulamaz ikiy\"uzl\"ul\"ukten\\
Can\i$\;$ay\i rmaya kalkarlar bedenden;\\
Horoz gibi tepemde testere olsa\\
Akl\i m\i n kafas\i n\i$\;$keser atar\i m ben. \\
$\; $\\
$\;\;\;\;\;\;\;\;\;\;\;\;\;\;\;\;\;\;\;\;\;\;\;\;\;\;\;\;\;\;\;\;\;\;\;\;\;\;\;\;\;\;\;\;\;\;\;\;$ \"O. Hayyam
\end{verse}
\end{table}

\end{document}